\documentclass[a4paper,12pt]{article}
\usepackage{amsfonts,amsmath}
\usepackage{amssymb}
\usepackage[mathscr]{eucal}
\usepackage{amsthm}
\usepackage{amscd}
\usepackage[T2A]{fontenc}
\usepackage[dvips]{color}
\usepackage{graphicx,epsfig}
\usepackage[cp1251]{inputenc}
\usepackage[T2A]{fontenc}

\newcommand{\df}{\stackrel{\textrm{def}}{=}}
\newcounter{rem}[section]

\def\text{\mbox}

\def\varkappa{\kappa}
\def\suml{\sum\limits}
\def\maxl{\max\limits}
\def\minl{\min\limits}
\def\intl{\int\limits}
\def\supl{\sup\limits}
\def\infl{\inf\limits}

\def\prodl{\prod\limits}

\def\dpfrac{\displaystyle\frac}

\def\text{\mbox}

\def\dzeta{\zeta}
\def\suml{\sum\limits}
\def\maxl{\max\limits}
\def\minl{\min\limits}
\def\intl{\int\limits}
\def\supl{\sup\limits}
\def\infl{\inf\limits}
\def\liminfl{\liminf\limits}

\def\prodl{\prod\limits}
\def\dip{\displaystyle}

\def\dpfrac{\dip\frac}
\title{Models with varying structure}
\author{Brodsky B.E., Darkhovsky B.S.}
\date{}

\begin{document}
\sloppy \maketitle

\abstract{{\bf \centerline {Abstract}}

\vspace{1cm}

In this paper the problems of the retrospective analysis of models
with time-varying structure are considered. These models include
contamination models with randomly switching parameters and
multivariate classification models with an arbitrary number of
classes. Our main task here is to classify observations with
different stochastic generation mechanisms. A new classification
method is proposed. We analyze its properties both theoretically
and empirically. The asymptotic optimality of the propodsed method
(by the order of convergence to zero of the estimation error) is
also established. At the end of the paper we consider multivariate
change-in-mean models and multivariate regression models.

}

\bigskip
{\bf 1. Introduction}

The previous papers of the authors (see, e.g.,  "Statistical analysis of models with varying structure" (Applied Econometrics, 2015, in
Russian), "Multivariate models with varying structure: a binary case" (Review of Applied and Industrial Mathematics, 2016, in Russian) were
devoted to the main particular cases of the general problem: how to split univariate  mixtures of probabilistical distributions and to perform
multivariate classification with only two classes of observations (ordinary observations and outliers). In this paper we consider the general
problem of multivariate classification with an arbitrary number of classes of observations.

First, let us mention previous important steps into this field.
Models with switching regimes have a long pre-history in statistics
(see, e.g., Lindgren (1978)). A simple switching model with two
regimes has the following form:
$$\begin{array}{ll}
& Y_t=X_t\beta_1+u_{1t} \quad \text{for the 1st regime  } \\
& Y_t=X_t \beta_2+u_{2t} \quad \text{for the 2nd regime }.
\end{array}
$$

For models with endogenous switchings usual estimation techniques
for regressions are not applicable. Goldfeld and Quandt (1973)
proposed {\it regression models with Markov switchings}. In these
models probabilities of sequential switchings are supposed to be
constant. Usually they are described by the matrix of probabilities
of switchings between different states.

Another modification of the regression models with Markov switchings
was proposed by Lee, Porter (1984). The following transition matrix
was studied:
$$
\Lambda=[p_{ij}]_{i,j=0,1}, \quad p_{ij}=P\{I_t=j|I_{t-1}=i\}.
$$

Lee and Porter (1984) consider an example with railway transport in
the US in 1880-1886s which were influenced by the cartel agreement.
The following regression model was considered:
$$
log P_t=\beta_0+\beta_1 X_t+\beta_2 I_t+u_t,
$$
where $I_t=0$ or $I_t=1$ in dependence of 'price wars' in the
concrete period.

Cosslett and Lee (1985) generalized the model of  Lee and Porter to
the case of serially correlated errors $u_t$.

Many economic time series occasionally exhibit dramatic breaks in
their behavior, associated with with events such as financial crises
(Jeanne and Mason, 2000; Cerra, 2005; Hamilton, 2005) or abrupt
changes in government policy (Hamilton, 1988; Sims and Zha, 2004;
Davig, 2004). Abrupt changes are also a prevalent feature of
financial data and empirics of asset prices (Ang and Bekaert, 2003;
Garcia, Luger, and Renault, 2003; Dai, Singleton, and Wei, 2003).

The functional form of the 'hidden Markov model' with switching
states can be written as follows:
$$
y_t=c_{s_t}+\phi y_{t-1}+\epsilon_t,
$$
where $s_t$ is a random variable which takes the values $s_t=1$ and
$s_t=2$ obeying a two-state Markov chain law:
$$
Pr(s_t=j|s_{t-1}=i,s_{t-2}=k,\dots,y_{t-1},y_{t-2},\dots)=Pr(s_t=j|s_{t-1}=i)=p_{ij}.
\eqno(ii)
$$

A model of this form with no autoregressive elements ($\phi=0$)
appears to have been first analyzed by Lindgren (1978) and Baum,
et al. (1980). Specifications that incorporate autoregressive
elements date back in the speech recognition literature to Poritz
(1982), Juang and Rabiner (1985), and Rabiner (1989).
Markov-switching regressions were first introduced in econometrics
by Goldfeld and Quandt (1973), the likelihood function for which
was first calculated by Cosslett and Lee (1985). General
characterizations of moment and stationarity conditions for
Markov-switching processes can be found in Tjostheim (1986), Yang
(2000), Timmermann (2000), and Francq and Zakoian (2001).

A useful review of modern approaches
 to estimation in Markov-switching models can be found in Hamilton (2005).

However, the mechanism of Markov chain modeling is far not unique in
statistical description of dependent observations. Besides Markov
models, we can mention martingale and copula approaches to dealing
with dependent data, as well as description of statistical
dependence via different coefficients of 'mixing'. All of these
approaches are interrelated and we must choose the most appropriate
method for the concrete problem. In this paper we choose the mixing
paradigm for description of statistical dependence.

Now let us mention some important problems which lead to stochastic models with switching regimes.

\bigskip {\it Splitting mixtures of probabilistic distributions}

In the simplest case we suppose that the d.f. of observations has
the following form:
$$
F(x)=(1-\epsilon)F_0(x)+\epsilon F_1(x),
$$
where $F_0(x)$ is the d.f. of ordinary observations; $F_1(x)$ is the
d.f. of abnormal observations; $0\le \epsilon <1$ is the probability
of obtaining an abnormal observation.

We need to test the hypothesis of statistical homogeneity (no
abnormal observations) of an obtained sample
$X^N=\{x_1,x_2,\dots,x_N\}$. If this hypothesis is rejected then
we need to classify this sample into sub-samples of ordinary and
abnormal observations.

\bigskip
{\it Estimation for regression models with abnormal observations}

The natural generalization of the previous model is the regression
model with abnormal observations
$$
Y=X\beta+\epsilon,
$$
where $Y$ is the $n\times 1$ vector of dependent observations; $X$
is the $n\times k$ matrix of predictors; $\beta$ is $k\times 1$
vector of regression coefficients; $\epsilon$ is the $n\times 1$
vector of random noises with the d.f. of the following type:
$$
f_{\epsilon}(x)=(1-\delta) f_0(x)+\delta f_1(x),
$$
where $0\le \delta <1$ is the probability to obtain an abnormal
observation; $f_0(x)$ is the density function of ordinary
observations; $f_1(x)$ is the density function of abnormal
observations. For example, in the model with Huber's contamination
[Huber, 1985]: $f_0(\cdot)={\cal N}(0,\sigma^2),\;f_1(\cdot)={\cal
N}(0,\Lambda^2)$, and $\Lambda>>\lambda>0$.

\bigskip
{\it Estimation for regression models with changing coefficients}

Regression models with changing coefficients is another
generalization of the contamination model. We suppose that
regression coefficients of this model can change (switch) from the
level $\beta_0$ to $\beta_1$ and the mechanism of this change is
random. We need to test the hypothesis about the absence of
switchings for each coefficient ($\epsilon=0$) and in the case of
rejection of this hypothesis to classify observations into
different groups.

We need again to test the hypothesis of statistical homogeneity of
an obtained sample and to divide this sample into sub-samples of
ordinary and abnormal observations if the homogeneity hypothesis is
rejected.

The goal of this paper is to propose methods which can solve these
problems effectively. Theoretically, we mean estimation of type 1
and type 2 errors in testing the statistical homogeneity hypothesis
and with estimation of contaminations parameters in the case of
rejectiong this hypothesis. Practically, we propose procedures for
implementation of these methods for univariate and multivariate
models.

Problems considered in this paper differ substantially from classical change-point problems in which we suppose that distances between various
regimes are big enough. In this paper we consider contamination models with coefficients changing in a random way.

The structure of this paper is as follows. In sections 2 and 3 we
consider univariate models with switching effects. In section 2
for binary mixtures of probabilistic distributions we prove
theorem 1 about exponential convergence to zero of type 1 error in
classification (to detect switches for a statistically homogenous
sample) as the sample size $N$ tends to infinity; theorem 2 about
exponential convergence to zero of type 2 error (vice versa, to
accept stationarity hypothesis for a sample with switches). In
section 3.3 we prove theorem 3 which establishes the lower bound
for the error of classification for binary mixtures. From theorems
2 and 3 we conclude that the proposed method is asymptotically
optimal by the order of convergence to zero of the classification
error.

Different generalizations of the proposed method for the case of
univariate models with multiple switching regimes and for
multivariate models with switching regimes are considered in
sections 3.4 and 3.5. Results of a detailed Monte Carlo study of
the proposed method for different stochastic models with switching
regimes are presented.

In section 4 we consider multivariate models. Multivariate
classification problems are considered in section 4.1. Section 4.2
deals with multivariate regression models.

\bigskip
 {\large \bf 2. Problems statement}

\bigskip
{\bf 2.1. Change-in-mean problems}

Suppose the d.f. of the observations is the binary mixture
$$
f(x)=(1-\epsilon)f_0(x)+\epsilon f_1(x),
$$
where the density functions $f_0(\cdot),\;f_1(\cdot)$ and the
parameter $\epsilon$ are unknown. We also suppose that
$$\mathbf
{E}_0(x)=\int x\,f_0(x)dx=0,\;\mathbf {E}_1(x)=\int x \,f_1(x)dx
=h\ne 0,
$$
where everywhere in this paper we denote by
$\mathbf{P}_0(\mathbf{E}_0)$ measure (mathematical expectation) of
the sequence $X^N$ under the condition $\epsilon=0$ (no 'abnormal'
observations.

The problem is to classify all obtained observations into
subsamples of ordinary data and outliers.

The estimation method is as follows:

1) From the initial sample $X^N$  compute the estimate of the mean
value:
$$
\theta_N=\dpfrac 1N \suml_{i=1}^N\,x_i
$$

2) Fix the numbers $0<\kappa<B$ and parameter $b\in \mathbb{B}\df
[\kappa, B]$ and classify observations as follows: if an observation
falls into the interval $(\theta_N-b,\theta_N+b)$, then we place it
into the sub-sample of ordinary observations, otherwise - to the
sub-sample of abnormal observations.

3) Then for each $b\in \mathbb{B}$ we obtain the following
decomposition of the sample $X^N$ into two sub-samples
$$\begin{array}{ll}
& X_1(b)=\{\tilde x_1,\tilde x_2,\dots,\tilde x_{N_1}\},\quad |\tilde x_i-\theta_N|< b, \\
& X_2(b)=\{\hat x_1,\hat x_2,\dots,\hat x_{N_2}\},\quad |\hat
x_i-\theta_N|\ge b
\end{array}
$$
Denote by $N_1=N_1(b),\,N_2=N_2(b),\,N=N_1+N_2$ the sizes of the
sub-samples $X_1$ and $X_2$, respectively.

The parameter  $b$  is chosen so that the sub-samples $X_1$ and
$X_2$ are separated in the best way. For this purpose, consider the
following statistic:
$$
\Psi_N(b)=\dpfrac 1{N^2}(N_2 \suml_{i=1}^{N_1}\,\tilde
x_i-N_1\suml_{i=1}^{N_2}\,\hat x_i).
$$

4) Define the boundary $C>0$ and compare it with the value $J=\max
|\Psi_N(b)|$ on the set $b\in \mathbb{B}$. If $J\le C$ then we
accept the hypothesis $H_0$ about the absence of abnormal
observations; if, however, $J>C$ then the hypothesis $H_0$ is
rejected.

Remark that testing the hypotheses $H_0,H_1$ does not require
knowledge of the distribution law of observations.

\bigskip
{\bf 2.2. Regression models with time-varying structure}

Here the following model of observations is considered:
$$
\mathbf Y=X\alpha+u_i=X(\mathbf\zeta\beta^0+\mathbf{(1-\zeta)}\beta^1)+\mathbf U,
$$
where

$\mathbf Y=(y_1,\dots,y_N)^{'}$  is a $N\times 1$ vector of dependent observations (here and below the sign $'$ denotes matrix transposition);

$X$  - $N\times k$ matrix of predictors;

$U$ - $N\times 1$ vector of centered random noises;

$\mathbf\alpha$  - $k\times 1$ vector of model coefficients,

$\mathbf\zeta$ - Bernoulli distributed random variable (independent from $\mathbf U$) with two states: $1$ with probability $(1-\epsilon)$ and
$0$ with probability $\epsilon$ for a certain unknown parameter $0<\epsilon<1$. Here $\beta^0\ne \beta^1$, $k$ is the number of model
coefficients.

In words we suppose that coefficients of this model can switch from the level $\beta^0$ into the level $\beta^1$, and the mechanism of these
switchings is random. We need to test the hypothesis of no switches in each coefficient ($\epsilon=0$).

Below we propose the method of solving this problem by means of
its reduction to the previous problem.

\vspace{0.5 cm} {\large \bf 3. Main results}

\bigskip
{\bf 3.1. Assumptions}

The results given below are based upon two main assumptions. The
first assumption is formulated in the form of a condition of
diminishing dependence between the past and the future of observed
processes as the distsnce between them increases. The second
condition takes the form of Cramer's assumption about the speed of
decrease of 'tails' of distributions.

A1.

a). {Mixing conditions}

On the probability space $(\Omega,\mathfrak {F},\mathbf{P})$ let
$\mathcal{H}_{1}$ and $\mathcal{H}_{2}$ be two $\sigma $-algebras
from $\mathfrak{F}$. Consider the following measure of dependence
between $\mathcal{H}_{1}$ and $\mathcal{H}_{2}$:
$$
 \psi (\mathcal{H}_{1}, \mathcal{H}_{2})  = \sup_{A \in
\mathcal{H}_{1},B \in \mathcal{H}_{2}, \mathbf{P}(A)
\mathbf{P}(B)\neq 0} \Big\vert \frac{\mathbf{P}(AB)}{\mathbf{P}(A)
\mathbf{P}(B)}-1 \Big\vert
$$

Suppose $\{y_n\},\,n\ge 1$ is a sequence of random variables defined
on ($\Omega,\mathfrak{F},\mathbf{P}$). Denote by
$\mathfrak{F}^{t}_{s}=\sigma \{y_{i}: s\le i\le t\}, 1\le s\le t<
\infty$ the minimal $\sigma $-algebra generated by random variables
$y_{i}, s \le i \le t$. Define
$$
 \psi (n)  = \sup_{t\ge 1} \psi (\mathfrak{F}^{t}_{1},
\mathfrak{F}^{\infty }_{t+n})
$$

We say that a random sequence $\{y_n\}$ satisfies the $\psi
$-\emph{mixing condition} if the function $\psi(n)$ (which is also
called the
 \emph{$\psi $-mixing coefficient}) tends to zero as $n$ goes to infinity.

The $\psi$-mixing condition is satisfied in most practical cases. In
particular, for a Markov chain (not necessarily stationary), if
$\psi(n)<1$ for a certain $n$, then $\psi(k)$ goes to zero at least
exponentially as $k\to\infty$ (see Bradley, 2005, theorem 3.3).

b) Nowadays, however, the notion of "weak dependence" of
observations is more often used:

\bigskip
{\bf Definition 2} (Doukhan, Louhichi, 1999). The sequence
$\{X_i\}$ is called $(\theta,{\cal L},\psi)$-weak dependent (or
simply $\psi$-weak dependent), if there exists a sequence
$\theta=(\theta_r)$ tending to zero as $r\to\infty$, and the
function $\psi$ with the argument $(f,h,n,m)\in {\cal L}_n \times
{\cal L}_m \times N^{2}$ such that for any sets of indices
$(i_1,\dots,i_n)$ and $(j_1,\dots,j_m)$  ($i_1\le \dots\le
i_n<i_n+r<j_1\le \dots\le j_m$):
$$
|Cov(f(X_{i_1},\dots,X_{i_n}),h(X_{j_1}\dots X_{j_m})| \le
\psi(f,h,n,m)\theta_r.
$$

It is often supposed that
$$
\theta_r=e^{-\beta r}, \qquad \beta>0.
$$

The 'weak dependence' condition is true in majority of practical
cases. In particular, Ango Nze, Doukhan (2004) showed that
$\psi$-weak dependence assumption generalizes conditions of mixing,
association , etc. , for Gaussian sequences and 'Bernoulli shifts'.
They proved that all ARMA and bilinear processes are $\psi$-weak
dependent. We can assume $\psi$- weak dependence while considering
all practically important cases in statistics.

A2.\textbf{Cramer condition}

We say that the sequence $\{y_n\}$ satisfies the \emph{uniform
Cramer condition} if there exists $T>0$ such that for each $i$,
$\mathbf{E}\exp(ty_i) < \infty$ for $|t|<T$.

For a centered sequence $\{y_n\}$ this condition is equivalent to
the following (see Petrov, 1987): there exist $g>0,\,H>0$ such that
$$
\mathbf{E}e^{ty_n}\le e^{\frac 12 g t^2} , \qquad |t|\le H,
$$
for all $n=1,2,\dots$.

\emph{We assume that conditions A1 and A2 hold true everywhere in
the paper.}

For any $x>0$ let us choose the number $\gamma(x)$ from the
following condition:
$$
\ln(1+\gamma(x))=\left \{
\begin{array}{ll}
& \dpfrac {x^2}{4g}, \qquad x\le gH \\
& \dpfrac {xH}4, \qquad x>gH,
\end{array}
\right.
$$
where $g,H$ are taken from the uniform Cramer condition.

For the chosen $\gamma(x)$, let us find such integer $\phi_0(x)\ge
1$ from the $\psi$-mixing condition that $\psi(l)\le\gamma(x)$ for
$l\ge \phi_0(x)$.

In the following theorem the exponential upper estimate for type 1
error  is obtained for the proposed method.

\bigskip
{\bf 3.2. Method}

Below we use the statistic $\Psi_N(b)$ defined in the previous
section. We note that it is a variant of the statistic that first
appeared in our papers and books (Brodsky, Darkhovsky, 1986, 1993,
2000) devoted to the analysis of change-point problems.
Methodologically, it ascends to Kolmogorov's test for detection
differently distributed random samples and to Hurst test in R/S
analysis.

\bigskip
{\bf Theorem 1.}

Let $\epsilon=0$. Suppose the d.f. $f_0(\cdot)$ is symmetric
w.r.t. zero and bounded. If $\psi$-mixing and Cramer's conditions
are satisfied then for any $0<\kappa<B$ there exists $C>0$ such
that the following estimate holds:
$$
\mathbf{P}_0\{\supl_{b\in\mathbb{B}}\,|\Psi_N(b)| >C\} \le L_1
\exp(-L_2(C)N),
$$
where the constants  $L_1,L_2>0$ do not depend on $N$.

However, if $\psi$-weak dependence and Cramer's conditions are satisfied then
$$
\mathbf{P}_0\{\supl_{b\in\mathbb{B}}\,|\Psi_N(b)| >C\} \le L_1
\exp(-L_2(C)\sqrt{N}),
$$
where, again, the constants  $L_1,L_2>0$ do not depend on $N$.

The proof of Theorem 1 is given in the Appendix.

\bigskip
Now consider characteristics of this method in case $\epsilon h\ne
0$. Here we again assume that $\mathbf{E}_0\,x_i=0, \,i=1,\dots,N$.

Put (for some fixed $\epsilon, h$)
$$
\begin{array}{ll}
&r(b)=\intl_{\epsilon h-b}^{\epsilon h+b}\,f(x)x dx,\quad d(b)=\intl_{\epsilon h-b}^{\epsilon h+b}\,f(x) dx\\[2mm]
&\Psi(b)=r(b)-\epsilon h d(b).
\end{array}
$$

In the following theorem type 2 error is studied.

\vspace{1cm} {\bf Theorem 2.}

1) Suppose $\psi$-mixing and Cramer's conditions are satisfied and
there exists $r^*=\supl_{b\in\mathbb{B}} r(b)$. Suppose also that
the density function $f(x)$ is continuous and there exists
$f^{''}(\cdot)\ne 0$. Then for $0< C<
\maxl_{b\in\mathbb{B}}\,|\Psi(b)|$ we have
$$
\mathbf{P}_{\epsilon}\{\maxl_{b\in\mathbb{B}}\,|\Psi_N(b)|\le
C\}\le L_1\exp(-L_2(\delta)N)).
$$
where $\delta=\maxl_{b\in\mathbb{B}}\,|\Psi(b)|-C
>0$..

2) If $\psi$-weak dependence condition is satisfied instead of
$\psi$-mixing, then
$$
\mathbf{P}_{\epsilon}\{\maxl_{b\in\mathbb{B}}\,|\Psi_N(b)|\le
C\}\le L_1\exp(-L_2(\delta)\sqrt{N}).
$$
where $\delta=\maxl_{b\in\mathbb{B}}\,|\Psi(b)|-C
>0$.

3) For solving estimation problems, we suppose that the underlying
model is
$$
f(x)=(1-\epsilon)f_0(x)+\epsilon f_0(x-h),    \eqno(!)
$$
where $0<\epsilon<1/2, h>0$ are unknown positive parameters.

Let us consider the equation:
$$
f(\epsilon h-b^*)=f(\epsilon h+b^*)   \eqno(!!)
$$

 Here we suppose that equation (!!)
has a unique root $b^{*}$ (for fixed $\epsilon, h$). Then
 $b_N^{*}\to b^{*}$ $\mathbf{P}_{\epsilon}$-a.s. as $N\to\infty$;,
where  $b_N^{*}>0$ is the estimate of $b^{*}$: $b_N^{*}\in
\text{arg}\maxl_{b\in\mathbb{B}}|\Psi_N(b)|$. Consider the
following estimates of $\epsilon$ and $h$:

$$\begin{array}{ll}
& \hat\epsilon_N \,\hat h_N=\theta_N \\
& \dpfrac {1-\hat\epsilon_N}{\hat\epsilon_N}=\dpfrac
{f_0(\theta_N-b_N^*-\hat h_N)-f_0(\theta_N+b_N^*-\hat
h_N)}{f_0(\theta_N+b_N^*)-f_0(\theta_N-b_N^*)}.
\end{array}
$$

Then the estimates $\hat\epsilon_N, \,\hat h_N$ converge
$\mathbf{P}_{\epsilon}$-a.s. to the true values of the parameters
$\epsilon, h$, respectively, as $N\to\infty$.

The proof of theorem 2 is given in the Appendix.

\bigskip
{\bf  Simulations}

We note that all constants in the above upper estimates of type 1 and type 2 errors are purely qualitative by their nature. Therefore
simulations of the proposed method are essential in the analysis of its properties.

In the first series of tests the following mixture model was
studied:
$$
f_{\epsilon}(x)=(1-\epsilon)f_0(x)+\epsilon f_0(x-h), \quad
f_0(\cdot)={\cal N}(0,1), \quad 0\le \epsilon <1/2.
$$

First, the critical thresholds of the decision statistic
$\max_{b\in\mathbb{B}} | \Psi_N(b)|$ were computed. For homogenous
samples of different size (i.e. without switches), p-quantiles of
the decision statistic were computed. For this purpose, a Gaussian
random sample with determined parameters was generated. After that
all steps of the above described method were done. The values of
the method's parameters: $\varkappa=0.04, B=50$.

The maximum of the absolute value of the decision statistic was
computed. This procedure was iterated 1000 times and the variation
series of the maximums of the absolute values of the decision
statistic was constructed. Then p-quantiles  (with $p=0.95$ and
$p=0.99$) in this series were computed. The obtained results are
given in Table 1.

\bigskip
{\bf Table 1. }

\bigskip
\hspace{-0.9cm }
\begin{tabular}{|c|c|c|c|c|c|c|c|c|c|c|}
\hline
$N$ & 50 & 100 &  300 & 500 & 800 & 1000 & 1200 & 1500 & 2000 \\
\hline
$\alpha=0.95$ & 0.1681 & 0.1213 & 0.0710 & 0.0534 & 0.044 & 0.0380 & 0.037 & 0.034 & 0.029 \\
\hline
$\alpha=0.99$ & 0.1833 & 0.1410 & 0.0869 & 0.0666 & 0.050 & 0.0471 & 0.0390 & 0.038 & 0.035 \\
\hline
\end{tabular}

\bigskip
In the second series of tests the quantile value for $p=0.95$ was
chosen as the critical threshold $C$ in experiments with
non-homogenous samples (for $\epsilon\ne 0$). For different sample
sizes in 1000 independent trials of each test, the estimate of
type 2 error $w_2$ (qi.e. the frequency of the event
$\maxl_{b\in\mathbb{B}} | \Psi_N(b)| <C$ for $\epsilon
>0$). The results are presented in table 2.

\bigskip
{\bf Table 2.}

\bigskip
\begin{tabular}{|c|c|c|c|c|c|c|c|c|}
\hline
$\epsilon=0.1 $ & \multicolumn{4}{|c|} {h=2.0} & \multicolumn{4}{|c|} {h=1.5} \\
\hline
$N$ & 300 & 500 &  800 & 1000 & 800 & 1200 & 2000 & 3000 \\
\hline
$C$ & 0.0710 & 0.0534 & 0.044 & 0.038 & 0.044 & 0.037 & 0.029 & 0.022 \\
\hline
$w_2$ & 0.26 & 0.15 & 0.05 & 0.02 & 0.62 & 0.42 & 0.16 & 0.03 \\
\hline
\end{tabular}

\vspace{1.5cm}
 {\bf 3.3. Asymptotic optimality}

\bigskip
Now consider the question about the asymptotic optimality of the
proposed method in the class of all estimates of the parameter
$\epsilon$. The a priori theoretical lower bound for the estimation
error of the parameter $\epsilon$ in the model with i.i.d.
observations with d.f. $f_{\epsilon}(x)=(1-\epsilon)f_0(x)+\epsilon
f_1(x)$ is given in the following theorem.

\bigskip
{\bf Theorem 3.} Let ${\cal M}_N$ be the class of all estimates of
the parameter $\epsilon$. Then for any $0< \delta < \epsilon$,
$$
\liminfl_{N\to\infty}\infl_{\hat\epsilon_N\in {\cal M}_N}\,\supl_{0<
\epsilon <1/2}\,\dpfrac 1N \ln
\mathbf{P}_{\epsilon}\{|\hat\epsilon_N-\epsilon|> \delta \} \ge
-\delta^2\,J(\epsilon),
$$
where $J(\epsilon)=\int\,[(f_0(x)-f_1(x))^2/ f_{\epsilon}(x)]\,dx$
is the generalized $\varkappa^2$ distance between densities $f_0(x)$
and $f_1(x)$ and $\mathbf{P}_{\epsilon}$ is the measure
corresponding to the density $f_{\epsilon}(x)$.

{\bf Proof.}

Remark that it suffices to consider consistent estimates of the
parameter $\epsilon$ (for non-consistent estimates the limit in the
left hand of the above inequality is equal to zero). This class is
not empty because of the method proposed in the paper.

Suppose $\hat\epsilon_N$ is any consistent estimate of $\epsilon$
and $0< \delta < \delta^{'}$. Consider the random variable
$\lambda_N=\lambda_N(x_1,\dots,x_N)=\mathbb{I}\{|\hat\epsilon_N-\epsilon|>
\delta \}$, where $\mathbb{I}(A)$ is the indicator of the set $A$.

Then for any $d>0$:
$$
\mathbf{P}_{\epsilon}\{|\hat\epsilon_N-\epsilon| > \delta
\}=E_{\epsilon}\lambda_N\ge \mathbf{E}_{\epsilon} (\lambda_N
\mathbb{I}\{f(X^N, \epsilon+\delta^{'})/ f(X^N,\epsilon) < e^d\}),
$$
where $f(X^N,\epsilon)$ is the likelihood function of the sample
$X^N$ of observations with the density function $f_{\epsilon}(x)$,
i.e.
$$
f(X^N,\epsilon)=\prodl_{i=1}^N\,[(1-\epsilon)f_0(x_i)+\epsilon
f_1(x_i)].
$$

Further,
$$\begin{array}{ll}
& \mathbf{E}_{\epsilon}(\lambda_N\mathbb{I}\{\dpfrac {f(X^N,\epsilon+\delta^{'})}{f(X^N,\epsilon)}<e^d\})\ge \\
& \ge e^{-d}\mathbf{E}_{\epsilon+\delta^{'}}(\lambda_N \mathbb{I}\{f(X^N,\epsilon+\delta^{'})/ f(X^N,\epsilon) < e^d \} \ge \\
& \ge
e^{-d}\,(\mathbf{P}_{\epsilon+\delta^{'}}\{|\hat\epsilon_N-\epsilon|>
\delta\}-\mathbf{P}_{\epsilon+\delta^{'}}\{f(X^N,\epsilon+\delta^{'})/
f(X^N,\epsilon)
> e^d\}).
\end{array}
$$

Since $\hat\epsilon_N$ is a consistent estimate,
$\mathbf{P}_{\epsilon+\delta^{'}}\{|\hat\epsilon_N-\epsilon| >
\delta \}\to 1$ as $N\to\infty$.

Let us consider the probability
$\mathbf{P}_{\epsilon+\delta^{'}}\{f(X^N,\epsilon+\delta^{'})/
f(X^N,\epsilon)
> e^d\}$. We have
$$\begin{array}{ll}
& \ln\dpfrac {f(X^N,\epsilon+\delta^{'})}{f(X^N,\epsilon)}=\suml_{i=1}^N\,\ln (1+\delta^{'}\dpfrac {f_1(x_i)-f_0(x_i)}{f_{\epsilon}(x_i)})= \\
& =\delta^{'}\suml_{i=1}^N\,\dpfrac
{f_1(x_i)-f_0(x_i)}{f_{\epsilon}(x_i)}+o(\delta^{'}).
\end{array}
$$
On the other hand,
$$
\mathbf{E}_{\epsilon+\delta^{'}}\dpfrac
{f_1(x_i)-f_0(x_i)}{f_{\epsilon}(x_i)}=\delta^{'}\int\,\dpfrac
{(f_1(x_i)-f_0(x_i))^2}{f_{\epsilon}(x_i)}dx_i=\delta ^{'}
J(\epsilon).
$$
Therefore, choosing
$d=N((\delta^{'})^2+\varkappa)J(\epsilon),\;\varkappa=o((\delta^{'})^2)$,
we obtain
$$
\mathbf{P}_{\epsilon+\delta^{'}}\{f(X^N,\epsilon+\delta^{'})/
f(X^N,\epsilon)
> e^d\}\to 0 \text {  as  } N\to\infty.
$$
Thus,
$$
\mathbf{P}_{\epsilon}\{|\hat\epsilon_N-\epsilon| > \delta \} \ge
(1-o(1))\,e^{-N\delta^2\,J(\epsilon)},
$$
or
$$
\liminfl_{N\to\infty}\infl_{\hat\epsilon_N\in {\cal M}_N}\,\supl_{0<
\epsilon <1/2}\,\dpfrac 1N \ln
\mathbf{P}_{\epsilon}\{|\hat\epsilon_N-\epsilon|> \delta \} \ge
-\delta^2\,J(\epsilon),
$$

Theorem 3 is proved.

Comparing results of theorems 2 and 3 we conclude that the proposed
method is asymptotically optimal by the order of convergence of the
estimates of a mixture parameters to their true values.

\vspace{1cm}
{\bf 3.4. Generalizations: non-symmetric distribution
functions}

Results obtained in theorems 1 and 2 can be generalized to the case
of non-symmetric distribution functions. Suppose the d.f.
$f_0(\cdot)$ is asymmetric w.r.t. zero. Then we can modify the
proposed method as follows.

1. From the initial sample $X^N=\{x_1,\dots,x_N\}$ compute the mean
value $\theta_N=\dpfrac 1N\,\suml_{i=1}^N\,x_i$ and the sample
$Y^N=\{y_1,\dots,y_N\};\;y_i=x_i-\theta_N$. Then we divide the
sample $Y^N$ into two sub-samples $I_1(b),\,I_2(b)$ as follows:
$$
y_i\in \left\{
\begin{array}{ll}
& I_1(b)=\{\tilde{y}_1,\dots,\tilde{y}_{N_1(b)}\},\quad -\phi(b)\le y_i\le b \\
& I_2(b)=\{\hat{y}_1,\dots,\hat{y}_{N_2(b)}\},\quad y_i>b \text{ or
} y_i<-\phi(b),
\end{array}
\right.
$$
where the function $\phi(b)$ is defined from the following
condition: $0=\intl_{-\phi(b)}^b\,y\,f_0(y)dy$,
$f_0(y)=f_0(x-\theta_N)$, $N=N_1(b)+N_2(b)$ and $N_1(b),N_2(b)$ are
sample sizes of $I_1(b),\,I_2(b)$, respectively.

2. As before we compute the statistic
$$
\Psi_N(b)=\dpfrac 1{N^2}(N_2(b) \suml_{i=1}^{N_1(b)}\,\tilde
y_i-N_1(b)\suml_{i=1}^{N_2(b)}\,\hat y_i).
$$

3. Then the value $J=\max_{b\in\mathbb{B}} |\Psi_N(b)|$ is compared
with the threshold $C$. If $J\le C$ then the hypothesis $H_0$ (no
abnormal observations) is accepted; if, however, $J>C$ then the
hypothesis $H_0$ is rejected and the estimate of the parameter
$\epsilon$ is constructed.

4. For this purpose, define the value $b_N^*$:
$$
b_N^* \in \arg\max_{b\in\mathbb{B}}|\Psi_N(b)|.
$$

Then
$$
\epsilon_N^*=N_2(b_N^*)/N.
$$

\vspace{0.5cm}

Consider application of this method for the study of the classic
$\epsilon$-contamination model:
$$
f_{\epsilon}(\cdot)=(1-\epsilon){\cal N}(\mu,\sigma^2)+\epsilon
{\cal N}(\mu, \Lambda^2), \quad \Lambda^2 >> \sigma^2,\quad
0\le\epsilon <1/2.
$$

For this model, the method described above has the form:

1. From the sample of observations $X^N=\{x_1,\dots,x_N\}$ the mean
value estimate $\hat\mu=\sum_{i=1}^N\,x_i/ N$ was computed.

2. The sequence $y_i=(x_i-\hat\mu)^2, \;i=1,\dots,N$ and its
empirical mean $\theta_N=\sum_{i=1}^N\,y_i/ N$ are computed.

3. Then for each $b\in\mathbb{B}$, the sample
$Y^N=\{y_1,\dots,y_N\}$ is divided into two sub-samples in the
following way: for $\theta_N (1-\phi(b))\le y_i\le \theta_N(1+b)$
put $\tilde y_i=y_i$ (the size of the sub-sample $N_1=N_1(b)$),
otherwise put $\hat y_i=y_i$ (the size of the sub-sample
$N_2=N_2(b)$). Here we choose the function $\phi(b)$ from the
following condition:
$$
\intl_{\theta_N (1-\phi(b))}^{\theta_N (1+b)}\,y f_0(y)dy=0,
$$
where $f_0(\cdot)=N(0,(1-\epsilon)^2\sigma^2)$.

From here we obtain:
$$
\phi(b)=1-\dpfrac b{e^b-1}.
$$

4. For any $b\in\mathbb{B}$, the following statistic is computed:
$$
\Psi_N(b)=\dpfrac 1{N^2}(N_2 \suml_{i=1}^{N_1}\,\tilde
y_i-N_1\suml_{i=1}^{N_2}\,\hat y_i).
$$
where $N=N_1+N_2,\,N_1=N_1(b),\,N_2=N_2(b)$ are sizes of sub-samples
of ordinary and abnormal observations, respectively.

5. Then, as above, the threshold $C>0$ is chosen and compared with
the value $J=\max_b |\Psi_N(b)|$. If $J\le C$ then the hypothesis
$H_0$ (no abnormal observations) is accepted; if, however, $J>C$
then the hypothesis $H_0$ is rejected and the estimate of the
parameter $\epsilon$ is constructed as follows.

Define the value $b_N^*$:
$$
b_N^* \in \arg\max_{b>0}|\Psi_N(b)|.
$$

Then
$$
\epsilon_N^*=N_2(b_N^*)/N.
$$

\textbf{Remark.} For estimation of the threshold, we use the
approach described in 2.1.3.

\bigskip In experiments the critical values of the statistic $\max_b | \Psi_N(b)|$ were computed. For this purpose, as above, for homogenous
samples (for $\epsilon=0$), $\alpha$-quantiles of the decision
statistic $\max_b | \Psi_N(b)|$ were computed
($\alpha=0.95,\,0.99$). The results obtained in 5000 trials of each
test are presented in table 3.

\bigskip
{\bf Table 3.}

\bigskip
\hspace{-0.8cm}
\begin{tabular}{|c|c|c|c|c|c|c|c|c|c|}
\hline
$N$ & 50 & 100 &  300 & 500 & 800 & 1000 & 1200 & 1500 & 2000 \\
\hline
$0.95$ & 0.3031 & 0.2330 & 0.1570 & 0.1419 & 0.1252 & 0.1244 & 0.1146 & 0.1107 & 0.1075 \\
\hline
$0.99$ & 0.3699 & 0.2862 & 0.1947 & 0.1543 & 0.1436 & 0.1331 & 0.1269 & 0.1190 & 0.1157 \\
\hline
\end{tabular}

\bigskip
The quantile value for $\alpha=0.95$ was chosen as the critical
threshold $C$ in experiments with non-homogenous samples (for
$\epsilon\ne 0$). For different sample sizes in 5000 independent
trials of each test, the estimate of type 2 error $w_2$ ( i.e. the
frequency of the event $\maxl_b | \Psi_N(b)| <C$ for $\epsilon
>0$) and the estimate $\hat\epsilon$ of the parameter $\epsilon$ were computed. The results are presented in tables 4 and 5.

\vspace{0.5cm} {\bf Table 4.}

\bigskip
\begin{tabular}{|c|c|c|c|c|}
\hline
$\Lambda=3.0 $ & \multicolumn{4}{|c|} {$\epsilon=0.05$}  \\
\hline
$N$ & 300 & 500 &  800 & 1000 \\
\hline
$C$ & 0.1570 & 0.1419 & 0.1252 & 0.1244  \\
\hline
$w_2$ & 0.27 & 0.15 & 0.06 & 0.04 \\
\hline
$\hat\epsilon$ & 0.064 & 0.056 & 0.052 & 0.05 \\
\hline
\end{tabular}

\bigskip
{\bf Table 5.}

\bigskip
\begin{tabular}{|c|c|c|c|c|c|}
\hline
$\Lambda=5.0 $ & \multicolumn{5}{|c|} {$\epsilon=0.01$}  \\
\hline
$N$ & 1000 & 1200 &  1500 & 2000 & 3000 \\
\hline
$C$ & 0.1244 & 0.1146 & 0.1107 & 0.1075 & 0.1019 \\
\hline
$w_2$ & 0.25 & 0.20 & 0.15 & 0.10 & 0.04 \\
\hline $\hat\epsilon$ & 0.0135 & 0.013 & 0.012 & 0.011 & 0.010
 \\
\hline
\end{tabular}

\vspace{1.5cm} {\bf 3.5. Generalizations: multiple switchings}

Suppose we obtain the data $X^N=\{x_1,\dots,x_N\}$, where the d.f.
of an observation $x_i$ can be written as follows:
$$
f(x_i)=(1-\epsilon_1-\dots-\epsilon_k)\,f_0(x_i)+\epsilon_1\,f_1(x_i)+\dots+\epsilon_k\,f_{k}(x_i),
$$
where $\epsilon_1\ge\epsilon_2\ge \dots \ge \epsilon_k \ge 0$, $0\le
\epsilon_1+\dots+\epsilon_k <1$, $|E_1f_1|< |E_2f_2|< \dots
<|E_{k}f_{k}|$.

In particular, we suppose the d.f. $f_0(x)$ is symmetric and
unimodal w.r.t. $E_0f_0$
$$
f(x_i)=(1-\epsilon_1-\dots-\epsilon_k)\,f_0(x_i)+\epsilon_1\,f_0(x_i-h_1)+\dots+\epsilon_k\,f_{0}(x_i-h_k),
$$
and $H=\epsilon_1 h_1+\dots+\epsilon_k h_k\ne 0$.

Our goal is to test the hypothesis $\epsilon_s=0,s=1,\dots,k$ (no
switches). In this section we denote by
$\mathbf{E}_i,\,i=0,1,\dots,k$, the mathematical expectation of
random variables corresponding to the d.f. with shift $\mathbf
{E}_if_i (\mathbf {E}_0f_0\df 0)$.

This model has the following sense. In the case of a binary
switching we have ordinary and abnormal observations. In the case of
multiple switchings abnormal observations are from different
classes. We do in analogy with the general form of this method.

1.1 From the initial sample $X^N$  compute the estimate of the mean
value:
$$
\theta_N=\dpfrac 1N \suml_{i=1}^N\,x_i
$$

1.2 Fix the numbers $0<\kappa<B$ and parameter $b\in \mathbb{B}\df
[\kappa, B]$ and classify observations as follows: if an observation
falls into the interval $(\theta_N-b,\theta_N+b)$, then we place it
into the sub-sample of ordinary observations, otherwise - to the
sub-sample of abnormal observations.

1.3. Then for each $b\in \mathbb{B}$ we obtain the following
decomposition of the sample $X^N$ into two sub-samples
$$\begin{array}{ll}
& X_1(b)=\{\tilde x_1,\tilde x_2,\dots,\tilde x_{N_1}\},\quad |\tilde x_i-\theta_N|< b, \\
& X_2(b)=\{\hat x_1,\hat x_2,\dots,\hat x_{N_2}\},\quad |\hat
x_i-\theta_N|\ge b
\end{array}
$$
Denote by $N_1=N_1(b),\,N_2=N_2(b),\,N=N_1+N_2$ the sizes of the
sub-samples $X_1$ and $X_2$, respectively.

The parameter  $b$  is chosen so that the sub-samples $X_1$ and
$X_2$ are separated in the best way. For this purpose, consider the
following statistic:
$$
\Psi_N(b)=\dpfrac 1{N^2}(N_2 \suml_{i=1}^{N_1}\,\tilde
x_i-N_1\suml_{i=1}^{N_2}\,\hat x_i).
$$

1.4. Define the boundary $C>0$ and compare it with the value $J=\max
|\Psi_N(b)|$ on the set $b\in \mathbb{B}$. If $J\le C$ then we
accept the hypothesis $H_0$ about the absence of abnormal
observations; if, however, $J>C$ then the hypothesis $H_0$ is
rejected.

2. As a result, we obtain two sun-samples: ordinary observations and
outliers at the first step of the algorithm.

3. Then we remove all found 'ordinary' observations from the sample
and repeat steps 1 and 2.

4. So we proceed further until a sub-sample without switches is obtained (i.e. the decision threshold $C$ is not exceeded).

Remark that the 1st type error for multiple switchings can be estimated like in the binary case (we do not formulate this result). As to the 2nd
type error (i.e. the probability that we stop at the 1st step of the method because the decision threshold is not exceeded) just observe that a
binary switch is a particular case of the general multiple switching situation (when all $\epsilon_i$ beginning from $i=2$ are equal to zero).

Therefore
$$\begin{array}{ll}
\mathbf{P}_{\epsilon}\{2nd \text{ type error, multiple switches}\}\le & \mathbf{P}_{\epsilon}\{2nd \text{ type error, binary case}\}\\
&\le L_1 \exp(-\beta(\delta,N)),
\end{array}
$$
for $0\le\delta\le \maxl_{b\in\mathbb{B}}\,|\Psi(b)|-C$.

\bigskip
{\bf Theorem 4.}

Suppose $0<C<\maxl_{b\in\mathbb{B}}\,|\Psi(b)|$. Then the 2nd type
error probability is estimated from above as follows:
$$
\mathbf{P}_{\epsilon}\{\text { 2nd type error }\}\le L_1 \exp(-\beta(\delta,N)),
$$
where $0\le\delta = \maxl_{b\in\mathbb{B}}\,|\Psi(b)|-C$

\bigskip{\bf Example}

Let us consider the following example. Suppose we have the model
with three classes of observations:
$$
f(x_i)=(1-\epsilon_1-\epsilon_2)\,f_0(x_i-h_1)+\epsilon_1\,f_0(x_i-h_2)+\epsilon_2\,f_0(x_i-h_3),
\qquad i=1,\dots,N,
$$
where $f_0(\cdot)={\cal N}(0,1)$; $x_i$ are i.r.v.'s.

The problem is to test the stationarity hypothesis: $H_0: \epsilon_1=\epsilon_2=0$.

Concretely, in this model the following parameters were chosen:
$$\begin{array}{ll}
& \epsilon_1=0.3;\;\epsilon_2=0.15 \\
& h_1=1,\;h_2=3,\;h_3=7.
\end{array}
$$

In experiments we estimated the type 2 error probability $\hat
w_2$.

The following results were obtained (each cell of this table is the
average in 1000 replications):

\bigskip
{\bf Table 6.}

\bigskip
\begin{tabular}{|c|c|c|c|c|c|c|c|}
  \hline
  $N$ & 100 & 200 & 300 & 500 & 700 & 1000 & 1500 \\
  \hline
  $\hat w_2$ & 0.116 & 0.090 & 0.070 & 0.048 & 0.036 & 0.016 & 0.010 \\
  \hline
\end{tabular}

\vspace{1.5cm} {\large \bf 4. Multivariate models}

\bigskip
{\bf 4.1. Multivariate classification}

\bigskip
{\bf Binary mixtures}

Now let us consider the multivariate classification problem with
binary mixtures. Suppose multivariate observations are of the
following type:
$$
{\mathbf{\cal Y}}^N=\{{\mathbf Y}^n\}_{n=1}^N,\,\,{\mathbf
Y}^n=(y_n^1,\dots,y_n^k).
$$

The multivariate density function of the vector $\mathbf Y^n$ is
$$
f(\mathbf Y^n)=(1-\epsilon)f_0(\mathbf Y^n)+\epsilon f_1(\mathbf
Y^n),
$$
where $f_0(\cdot),\,f_1(\cdot)$ are the d.f.'s of ordinary and
abnormal observations, respectively; the d.f. $f_0(\cdot)$ is
supposed to be symmetric w.r.t. its mean vector.

First, let us consider the case $\mathbf{E}_1(\mathbf
Y^n)=\mathbf{a}\ne 0$, i.e. changes in mean of abnormal
observations. Remark that the baseline "change-in-mean" problem is
usually considered in many methods of 'cluster analysis' in which
different distances between multivariate 'points' of
characteristics (even without references to density functions and
mathematical expectations of observations) are considered.

The method can be formulated in analogy with the univariate case:

1) From the initial sample ${\mathbf{\cal Y}}^N$  compute the
estimate of the mean value:
$$
\theta_N=\dpfrac 1N \suml_{i=1}^N\,\mathbf Y^i.
$$
2) Fix the parameter $b>0$ and classify observations as follows:

if $\|\mathbf Y^i-\theta_N\| \le b$, then we place $\mathbf Y^i$
into the sub-sample of ordinary observations $\{\tilde {\mathbf
Y}^i\}$;

if $\|\mathbf Y^i-\theta_N\| > b$, then we place $\mathbf Y^i$
into the sub-sample of abnormal observations $\{\hat {\mathbf
Y}^i\}$.

As a result, for each $b>0$ we obtain the decomposition of the
sample ${\mathbf{\cal Y}}^N$ into sub-samples of ordinary and
abnormal observations. Suppose the size of ordinary sub-sample is
$N_1(b)$ and the size of abnormal sub-sample is $N_2(b)$.

3) The parameter  $b$  can be chosen in order to separate the
sub-samples of ordinary and abnormal observations
($\{\mathbf{\tilde{Y}}^i\}$ and $\{\mathbf{\hat{Y}}^i\}$,
respectively) in the best way. For this purpose, consider the
following statistic:
$$
{\mathbf\Psi_N}(b)=\dpfrac 1{N^2}(N_2
\suml_{i=1}^{N_1}\,{\mathbf{\tilde
Y}}^i-N_1\suml_{i=1}^{N_2}\,{\mathbf{\hat Y}}^i).
$$

4) Define the boundary $C>0$ and compare it with the value
$J=\maxl_{b\in\mathbb{B}} \|\Psi_N(b)\|$. If $J\le C$ then we accept
the hypothesis $H_0$ about the absence of abnormal observations; if,
however, $J>C$ then the hypothesis $H_0$ is rejected.

For this method, in analogy with the univariate case we can
formulate results about type 1 and type 2 eroror probabilities.
For example, the exponential upper estimate for type 1 error
probability is formulated as follows:

Let $\epsilon=0$. Suppose the d.f $f_0(\cdot)$ is symmetric w.r.t.
zero and bounded. Then for all $\kappa, B:0<\kappa<B$ there exists
$C>9$ such that
$$
\mathbf{P}_0\{\supl_{b\in \mathbb{B}}\,\|\mathbf {\Psi}_N(b)\|\le
C\}\le L_1\exp(-\beta(C,N)),
$$
where $\mathbb{B}=[\kappa, B]$.

Imitation modeling

In this example the following multivariate Gaussian model was considered:
$$
f(\mathbf X)=(1-\epsilon)f_0(\mathbf X)+{\epsilon}f_1(\mathbf X).
$$
where $f_0(\mathbf X)$ is the two-dimensional Gaussian d.f. with the vector of means $\mu_1=(0\; 0)^{'}$ and the covariance matrix
$Cov(x_i)=\left(
\begin{array}{cc}
  0.745 & -0.07 \\
  -0.07 & 0.01 \\
\end{array}
\right) $, and $f_1(\mathbf X)$ is the two-dimensional Gaussian d.f. with the vector of means $\mu_2=(0\; 0.25)^{'}$ and the same correlation
matrix. Here $\epsilon=0.2$.

In this model it is a priori known that switchings occur in the second coordinate of observations. Therefore from the beginning we consider this
second coordinate (which is connected with the first coordinate in virtue of out two-dimensional model).

First, the critical thresholds of the decision statistic $\max_{b\in\mathbb{B}} | \Psi_N(b)|$ were computed. For homogenous samples of different
size (i.e. without switches), p-quantiles of the decision statistic were computed. For this purpose, a Gaussian random sample with determined
parameters was generated. After that all steps of the above described method were done. The values of the method's parameters: $\varkappa=0.04,
B=50$.

The maximum of the absolute value of the decision statistic was computed. This procedure was iterated 1000 times and the variation series of the
maximums of the absolute values of the decision statistic was constructed. Then p-quantiles  (with $p=0.95$ and $p=0.99$) in this series were
computed. The obtained results are given in Table 7.

\bigskip
{\bf Table 7. }

\bigskip
\hspace{-0.9cm }
\begin{tabular}{|c|c|c|c|c|c|c|c|c|c|}
\hline
$N$ & 50 & 100 &  200 & 300 & 500 & 700 & 1000 & 1500 \\
\hline
$\alpha=0.95$ & 0.0066 & 0.0059 & 0.0041 & 0.0037 & 0.0027 & 0.0024 & 0.0019 & 0.0016 \\
\hline
$\alpha=0.99$ & 0.014 & 0.0083 & 0.0057 & 0.0045 & 0.0037 & 0.0036 & 0.0024 & 0.0020 \\
\hline
\end{tabular}

\bigskip
In the second series of tests the quantile value for $p=0.95$ was
chosen as the critical threshold $C$ in experiments with
non-homogenous samples (for $\epsilon\ne 0$). For different sample
sizes in 1000 independent trials of each test, the estimate of
type 2 error $w_2$ (i.e. the frequency of the event
$\maxl_{b\in\mathbb{B}} | \Psi_N(b)| <C$ for $\epsilon
>0$). The results are presented in table 8.

\bigskip
{\bf Table 8.}

\bigskip
\begin{tabular}{|c|c|c|c|c|c|c|c|}
\hline
$N$ & 100 & 200 &  300 & 500 & 700 & 1000 & 1500 \\
\hline
$C$ & 0.0059 & 0.0041 & 0.0037 & 0.0027 & 0.0024 & 0.0019 & 0.0016 \\
\hline
$w_2$ & 0.110 & 0.019 & 0.002 & 0 & 0 & 0 & 0 \\
\hline
\end{tabular}

\bigskip
 The results obtained witness about the fact that the quality of this method increases with the growing sample size. Here: $w2$ is the
frequency of type 2 error, $C$ is the decision threshold.

\bigskip In analogy with the univariate case we can generalize this method to the case of multiple switchings.

\vspace{1cm} {\bf 4.2. Switching regressions}

Let us first remind the considered model of observations:
$$
\mathbf Y=X\mathbf\mathbf\beta+\mathbf
U=X(\mathbf\dzeta\beta_0+\mathbf{(1-\dzeta)}\beta_1)+\mathbf U,
$$
where

$\mathbf Y$  is a $N\times 1$ vector of dependent observations $y_1,y_2,\dots,y_N$;

$X$  is a $N\times k$ matrix of predictors;

$\mathbf U$  is a $N\times 1$ vector of centered random noises $u_1,u_2,\dots,u_N$;

$\mathbf\dzeta$ is a $k\times 1$ vector of r.v.'s $\dzeta_1,\dzeta_2,\dots, \dzeta_k$ independent of $u_1,\dots,u_N$ and identically distributed
according to Bernoulli law:
$$
\mathbf P\{\dzeta_j=1\}=1-\mathbf P\{\dzeta_j=0\}=\epsilon_j,
\quad j=1,2,\dots,k,
$$
for certain unknown parameters $0<\epsilon_j<1,\;j=1,\dots,k$,

$\mathbf {1}$ - $k\times 1$ vector composed of 1's.

Here $\beta_0\ne\beta_1,\; k$ - dimensionality of the vector of coefficients $\mathbf \alpha$ of the model.

For solving this problem, consider the OLS estimate of the vector
$\beta$ (here and below $'$ is the symbol of transposition):
$$
\hat\beta=(X'X)^{-1}X'\mathbf Y=\mathbf\dzeta\beta_0+\mathbf{(1-\dzeta)}\beta_1+(X'X)^{-1}X'\mathbf U.
$$

Since the sequence of noises $\mathbf U$ is centered, the problem is reduced to the above considered problem of detection of switches in the
mean of an observed random vector. The matrix of predictors $X$ influences only the random component.

Formally, we need to introduce the following vector $\mathbf I=(1,1,\dots,1)$ ($N$ units) and consider
$$
\mathbf{\tilde\beta}=[\dzeta\beta_0+(1-\dzeta)\beta_1]\,I+(X^{'}X)^{-1}X^{'}\mathbf U\,\mathbf I.
$$

Then the $(k\times N)$ matrix $\tilde \beta$ consists of $N$ columns of $k\times 1$ vectors with  means $\beta_0$ and $\beta_1$ changing in a
random manner. Each component $j=1,\dots,k$ of these vectors $\tilde \beta_i^j,\;i=1,\dots,N$ is therefore a univariate random sequence
$$
\tilde\beta_i^j=[\dzeta\beta_0^j+(1-\dzeta)\beta_1^j]_i+\xi_i^j, \qquad i=1,\dots,N,
$$
where
$$
\xi_i^j=((X^{'}X)^{-1}X^{'}\mathbf U\,\mathbf I)_i^j.
$$

So the problem of detection of changes in regression coefficients is reduced to the above considered problem of detection switches in the mean
value of a univariate random sequence. Remark that the uniform Cramer and the $\psi$-mixing conditions are still satisfied for the process
$\xi_i^j,\,i=1,\dots,N$. As $\mathbf{E}u_i\equiv 0$ we get that there exist constants $g_1>0,\,H_1>0$ such that
$$
Ee^{t\,\xi_i^j}\le e^{\dpfrac 12 g_1t^2},\quad |t|\le H_1,
$$
for all $i=1,\dots,N,\,j=1,\dots,k$. Moreover, we choose the number $m_0(\cdot)$ from the $\psi$-mixing condition for $\xi_i^j,\,i=1,\dots,N$:
for any chosen number $\gamma(x)>0$: $\psi(l)\le \gamma(x)$ for $l\ge m_0(x)$.

For testing the hypothesis of no switches we again consider the decision statistic $\Psi_N(b)$ and compare the maximum of its module with the
decision threshold $C>0$. Then the following theorem holds:

Formally, we can write
$$
\tilde\alpha=[\mathbf\dzeta\beta_0+\mathbf{(1-\dzeta)}\beta_1]+(X^{'}X)^{-1}X^{'}\mathbf U.
$$

Each component $j=1,\dots,k$ of these vectors $\tilde \alpha_i^j,\;i=1,\dots,N$ is therefore a univariate random sequence
$$
\tilde\alpha_i^j=[\dzeta_i^j\beta_0+(1-\dzeta_i^j)\beta_1]_i+\xi_i^j, \qquad i=1,\dots,N,
$$
where
$$
\xi_i^j=((X^{'}X)^{-1}X^{'}\mathbf U\,\mathbf I)_i^j.
$$

So the problem of detection of changes in regression coefficients is reduced to the above considered problem of detection switches in the mean
value of a univariate random sequence. Remark that the uniform Cramer and the $\psi$-mixing conditions are still satisfied for the process
$\xi_i^j,\,i=1,\dots,N$. As $\mathbf{E}\mathbf U\equiv 0$ we get that there exist constants $g_1>0,\,H_1>0$ such that
$$
Ee^{t\,\xi_i^j}\le e^{\dpfrac 12 g_1t^2},\quad |t|\le H_1,
$$
for all $i=1,\dots,N,\,j=1,\dots,k$. Moreover, we choose the number
$m_0(\cdot)$ from the $\psi$-mixing condition for
$\xi_i^j,\,i=1,\dots,N$: for any chosen number $\gamma(x)>0$:
$\psi(l)\le \gamma(x)$ for $l\ge m_0(x)$.

For testing the hypothesis of no switches we again consider the
decision statistic $\Psi_N(b)$ and compare the maximum of its module
with the decision threshold $C>0$. Then the following theorem holds:

\bigskip
{\bf Theorem 5.}

Suppose $\epsilon=0$, the d.f. of each component of the vector $\mathbf U$ is symmetric w.r.t. zero and the $\psi$-mixing and the uniform Cramer
conditions for $\xi_i^j,\,i=1,\dots,N$ are satisfied. Then for any threshold $C>0$ the following upper estimate for the 1st type error
probability holds:
$$
\mathbf{P}_0\{\maxl_{b\in\mathbb{B}}\,|\Psi_N(b)| >C\} \le L_1 \exp(-\beta(C,N)),
$$
where the function $\beta(C,N)$ is defined in the proof of Theorem 1.

The proof of theorem 5 is based upon the same ideas as the proof of theorem 1. Therefore it is omitted here.

\bigskip
{\bf 5. Simulations}

{\bf To 4.2}

In the following example the regression model with one deterministic
predictor was considered:
$$
y_i=c_1+c_2*i+u_i, \quad u_i\sim N(0;1),\quad i=1,\dots,n.
$$
$$
\xi\sim U[0;1]
$$
$$
\gamma=[c_1;\;c_2]=\left\{
\begin{array}{ll}
& \beta_1,\quad \epsilon<\xi\le 1 \\
& \beta_2,\quad 0\le\xi\le\epsilon
\end{array}
\right.
$$

\vspace{2cm} {\bf Table 9.}

\bigskip
\begin{tabular}{|c|c|c|c|c|}
\hline
$\epsilon=0.05 $ & \multicolumn{4}{|c|} {$\beta_1=[1; 1],\;\beta_2=[1; 2]$}  \\
\hline
$N$ & 300 & 500 &  800 & 1000 \\
\hline
$C$ & 0.07 & 0.05 & 0.04 & 0.03  \\
\hline
$w_2$ & 0.87 & 0.59 & 0.14 & 0.004 \\
\hline
$\hat\epsilon$ & 0.08 & 0.059 & 0.052 & 0.05 \\
\hline
\end{tabular}

\bigskip
{\bf Table 10.}

\bigskip
\begin{tabular}{|c|c|c|c|c|}
\hline
$\epsilon=0.1 $ & \multicolumn{4}{|c|} {$\beta_1=[1; 1],\;\beta_2=[1; 1.5]$}  \\
\hline
$N$ & 300 & 500 &  800 & 1000 \\
\hline
$C$ & 0.07 & 0.05 & 0.04 & 0.03  \\
\hline
$w_2$ & 0.83 & 0.65 & 0.13 & 0.0 \\
\hline
$\hat\epsilon$ & 0.15 & 0.12 & 0.102 & 0.10 \\
\hline
\end{tabular}

\bigskip
 {\large \bf  Conclusion}

In this paper we considered the problems of the retrospective analysis of models with time-varying structure. These models include contamination
models with randomly switching parameters and multivariate classification models with an arbitrary number of classes. Our main task here is to
classify observations with different stochastic generation mechanisms. We propose a new classification method and analyze its properties both
theoretically and empirically. It was proved that type 1 and type 2 errors of the proposed method converge to zero exponentially as the sample
size $N$ tends to infinity. The asymptotic optimality of the proposed method follows from theorem 3. In this theorem the theoretical lower bound
for the error of estimation of the model's parameters was established. This bound is attained for the proposed method (by the order of
convergence to zero of the estimation error). Then we consider generalizations of the proposed method to the case of non-symmetric d.f.'s of
ordinary observations and to the case of an arbitrary number of classes of observations with different stochastic generation mechanisms. The
multivariate models with time-varying structure are considered at the end of this paper. Here we consider multivariate change-in-mean models and
multivariate models with time-varying regression coefficients.

\bigskip
{\bf Appendix. Proofs of theorems}

{\bf Theorem 1. Proof.}

In subsequent considerations we use many times the following inequality which will be proved first: let $S_n=\suml_{k=1}^n \xi_k$, where
$\{\xi_k\}_{k=1}^{\infty}$ is the sequence of r.v.'s, satisfying conditions A1 and A2 (whether {$\psi$-mixing or $\psi$-weak dependence), and
$\mathbf{E}\xi_k\equiv 0$.

Then under A1(a) and A2, for sufficiently large  $N$, the following inequality holds true:
$$
\begin{array}{ll}
&\mathbf{P}\left\{|S_{\scriptscriptstyle N}|/N>x\right\}\le A(x)\exp\left(-B(x)\gamma N\right)
\end{array}
\eqno(*)
$$
where positive functions $A(\cdot),\,B(\cdot)$  can be computed explicitly.

Under A1(b) and A2, however,
$$
\begin{array}{ll}
&\mathbf{P}\left\{|S_{\scriptscriptstyle N}|/N>x\right\}\le A(x)\exp\left(-B(x)\gamma \sqrt{N}\right),
\end{array}
\eqno(**)
$$
where,again, positive functions $A(\cdot),\,B(\cdot)$  can be computed explicitly.

In the sequel we introduce the following notation:
$$
\beta(x,N)=\{
\begin{array}{ll}
& B(x)\gamma N, \quad\text{  in case of $\psi$-mixing}\\
& B(x)\gamma\sqrt{N}, \quad \text{ in case of $\psi$-weak dependence}
\end{array}.
$$
Then
$$
\mathbf{P}\left\{|S_{\scriptscriptstyle N}|/N>x\right\}\le A(x)\exp\left(-\beta(x,N))\right)   \eqno(***)
$$

Under $\psi$-mixing and Cramer's conditions this inequality was proved in Brodsky, Darkhovsky (2000). Here we prove it under assumptions A1(b)
and A2.

Under "weak dependence" assumption we proceed from Roussas-Ionnides inequality (see, e.g., Hwang, Shin (2014)) in the following form. Let
$\gamma$ be a certain large number: $\gamma\ge 2$. Suppose $(p_1,\dots,p_{\gamma})$ and $q_{\gamma}$ are positive numbers such that
$$
\dpfrac 1{p_1}+\dots+\dpfrac 1{p_{\gamma}}=\dpfrac 1{q_{\gamma}}<1,
$$

We assume that $\xi_1,\dots,\xi\-{\gamma}$ is a weak dependent sequence with the function $\theta_r$ in Definition 2. Suppose that
$$
E\|\xi_i\|^{p_i}< \infty, \qquad p_i>1,
$$
for $i=1,\dots,\gamma-1$, $\gamma\ge 2$. The analogous boundedness conditions are imposed on the functions $h$ and $f$ and their first
derivatives.

Then the following inequality holds:
$$
|E|\prodl_{i=1}^{\gamma}\,\xi_i|-\prodl_{i=1}^{\gamma}\,E[\xi_i]|\le B(\gamma-1)\,\theta_r^{1-\frac
1{q_{\gamma}}}\,\prodl_{i=1}^{\gamma}\,\|\xi_i\|_{p_i},
$$
where the constant $B$ does not depend on $\gamma$ and $r$.

For the proof of theorem 1 we split the sum $S_n$ into $\phi$ terms of the following types (choice of $\phi$ is explained below):
$$
S_n=S_n^1+\dots+S_n^{\phi},
$$
where
$$
S_n^i=\xi(i)+\xi(i+\phi(x))+\dots
 +\xi(i+\phi(x)[\frac {n-i}{\phi(x)}]),
$$
и $i=1,2,\dots,\phi(x)$.

Then
$$
P\{|S_n|/n\ge x\}\le \phi(x)\,\maxl_{1\le i\le \phi}\,P\{|S_n^i\ge (k(i)-1)x\}.
$$

We need to obtain the exponential upper estimate for the probability $P\{Z_k>x\}$, where
$$
Z_k=\suml_{j=1}^k\,\xi(i+\phi j).
$$

From Chebyshev's inequality we obtain
$$
P\{Z_k>x\}\le e^{-tx}Ee^{tZ_k}.
$$

From the Roussas-Ionnides inequality we have:
$$\begin{array}{ll}
&P\{|S_n^i|\ge C\dpfrac {Cn}{m_N(C)}\}\le \left \{
\begin{array}{ll}
&\exp(-\dpfrac {nC^2}{2gm_N^2(C)}),\quad C\le gT \\
& \exp(-\dpfrac {nCT}{2m_N(C)}), \quad C>gT
\end{array}
\right.\\
& +B(n-1)\theta^{1-\dpfrac 1{q_n}}(m_N(C))\left\{
\begin{array}{ll}
& \exp(-\dpfrac {C^2n^2}{2g\sum_{j=1}^np_j}),\quad C\le gT\dpfrac
{\sum_1^np_j}n \\
&\exp(-\dpfrac{nTC}2)\quad C> gT\dpfrac {\sum_1^np_j}n.
\end{array}
\right.
\end{array}
$$

Consider the second term in the right hand:
$$
B(n-1)\theta^{1-\dpfrac 1{q_n}}(m_N(C))\left\{
\begin{array}{ll}
& \exp(-\dpfrac {C^2n^2}{2g\sum_{j=1}^np_j}),\quad C\le gT\dpfrac
{\sum_1^np_j}n \\
&\exp(-\dpfrac{nTC}2)\quad C> gT\dpfrac {\sum_1^np_j}n.
\end{array}
\right.
$$

The direct calculation of the munimum of the function $ \theta^{1-\dpfrac 1{q_n}} \exp(-\dpfrac {C^2n^2}{2g\sum_{j=1}^np_j})$ , dependent on the
arguments $p_1,\dots,p_{n-1},q_n$, on condition that
$$
\dpfrac 1{p_1}+\dots+\dpfrac 1{p_{n-1}}+\dpfrac 1{q_n}=1
$$
gives
$$
q_n^*\sim n^2,\,p_i\sim n, i=1,\dots,n-1.
$$

Therefore $1-\dpfrac 1{q_N^*}>1/2$ for large enough  $n$. This fact yields the estimate
$$
P\{|S_n^i|/n\ge C\}\le m_N(C)\exp(-\dpfrac n{m_N(C)}C)+B(n-1)\exp(-(m_N(C))\beta/2), \qquad \beta>0.
$$

Then we choose $m_N(C)\sim \sqrt{N}$ and obtain
$$
P\{|S_n^i|/n\ge C\}\le L_1\exp(-L_2(C)\sqrt{N}).
$$

So in the case of $\psi$-weakly dependent variables we need to choose $m_N(C)\sim \sqrt{N}$.

The function $\theta_{\phi}$ exponentially converges to zero with the increase of $\phi$ (this fact holds true in most cases):
$$
\theta_{\phi}\le e^{-\beta \phi}, \quad \beta>0,
$$

Then we can take $\phi\sim k$ and conclude that
$$ P\{|Z_k|/k>x\}\le
L_1\exp(-L_2(x)\sqrt{k}),
$$
where $L_1,L_2$ are constants not dependent on $k$.

Therefore,
$$
P\{S_n/n\}\le L_1\exp(-L_2(x)\sqrt{n}),
$$
and
$$
\alpha_N= P_0\{\maxl_{[\alpha N]\le l\le N}\,\|Y_N(l)\| > C \} \le L_1\exp(-L_2(C)\sqrt{N}).
$$
where as before $L_1,L_2$ are constants not dependent on $N$.

For  $\Psi_N(b)$ we can write:
$$
\Psi_N(b)=\left(N\,\suml_{i=1}^{N_1(b)}\,\tilde x^i-N_1(b)\suml_{i\in\mathcal{N}}\,x^i\right)/ N^2.   \eqno(1)
$$

Then
$$\begin{array}{ll}
\mathbf{P}_0\{\supl_{b\in\mathbf{B}}|\Psi_N(b)|>C\}&\le \mathbf{P}_0\{\supl_{b\in\mathbf{B}}|\suml_{i=1}^{N_1(b)}\,\tilde x^i| >\dpfrac
{CN}2\}\\[3mm]
&+\mathbf{P}_0\{|\suml_{i\in\mathcal{N}}\,x^i|>\dpfrac {CN}2\}.
\end{array}
\eqno(2)
$$

Further,
$$
\mathbf{P}_0\left\{\supl_{b\in\mathbf{B}}| \suml_{i=1}^{N_1(b)}\,\tilde x^i|
>\dpfrac C2 N\right\}\le\suml_{n=1}^N\,\mathbf{P}_0\left\{\supl_{b\in\mathbf{B}}\{|\suml_{i=1}^{n}\,\tilde x^i|
>\dpfrac C2 n\}\cap \{N_1(b)=n\}\right\}.
\eqno(3)
$$

Consider the function
$$
\Delta(b)=\intl_{|x|\le b}\,f_0(x)dx.
$$

the function $\Delta(b)$ is continuous and $\minl_{b\in\mathbf{B}}\Delta(b)\ge \intl_{|x|\le\kappa}f_0(x)dx\df u$.

Now let us split the segment $\mathbf{B}=[\kappa, B]$ into equal parts with the interval such that $|\Delta(b_i)-\Delta(b_{i+1})|\le u/2$. In
virtue of uniform continuity of $\Delta(b)$ such split is possible (here $\{b_i\}$ are bounds of this split).

Denote the number of subsegments by $R$, and subsegments themselvelvs by $\mathbf{B}_s,\,s=1,\dots,R$. Then
$$
\begin{array}{ll}
&\mathbf{P}_0\left\{\supl_{b\in\mathbf{B}}\{|\suml_{i=1}^{n}\,\tilde x^i|
>\dpfrac C2 n\}\cap \{N_1(b)=n\}\right\}\\
&=\mathbf{P}_0\left\{\maxl_s\supl_{b\in\mathbf{B}_s}\{|\suml_{i=1}^{n}\,\tilde x^i|
>\dpfrac C2 n\}\cap \{N_1(b)=n\}\right\} \\
&\le R\maxl_s\mathbf{P}_0\left\{\supl_{b\in\mathbf{B}_s}\{|\suml_{i=1}^{n}\,\tilde x^i|
>\dpfrac C2 n\}\cap \{N_1(b)=n\}\right\}
\end{array}
\eqno(4)
$$

Consider the fixed subsegment $\mathbf{B}_{i}\df [b_i,b_{i+1}]$. From definition of the numbers $N_1(b)$ we obtain for each
$b\in\mathbf{B}_{i}$:
$$
 N_1(b_i)/N\le N_1(b)/N \le N_1(b_{i+1})/N
$$

Now let us construct estimates for probabilities of the following events:
$$
|N_1(b_i)/N-\Delta(b_i)|\le u/4,\quad|N_1(b_{i+1})/N-\Delta(b_{i+1})|\le u/4
$$

First, let us estimate the probability of deviation of the r.v. $\theta_N$ from its mathematical expectation $\mathbf{E}_m\theta_N\equiv 0$.

For each $\gamma>0$ and sufficiently large  $N$, from inequality (***) it follows that
$$
\mathbf{P}_0\{|\theta_N|>\gamma\}\le A(\gamma)\exp\left(-B(\gamma)N\right) \eqno(5)
$$
in case A1(a) (mixing) and
$$
\mathbf{P}_0\{|\theta_N|>\gamma\}\le A(\gamma)\exp\left(-B(\gamma)\sqrt{N}\right) \eqno(5)
$$
in case a1(b)($\psi$-weak dependence).

From definition, $N_1(b)=\suml_{k\in\mathcal{N}}\mathbf{I}(|x^k-\theta_N |\le b)$.

Then for every fixed $r>0$ we obtain
$$
\mathbf{P}_0\{|x^k-\theta_N|\le b\}\le\mathbf{P}_0\{|x^k|\le b+r\}+\mathbf{P}_0\{|\theta_N|>r\} \eqno(6)
$$

Moreover,
$$
\mathbf{P}_0\{|x^k-\theta_N|\le b\}\ge \mathbf{P}_0\{|x^k|\le b-r\}-\mathbf{P}_0\{|\theta_N|>r\} \eqno(7)
$$

For each point $b_i$ from the split of the segment $\mathbb{B}$ we obtain for sufficiently large $N$ (see Brodsky, Darkhovsky (2000)):
$$
\mathbf{P}_0\{|\frac{1}{N}\suml_{k\in\mathcal{N}}\left(\mathbb{I}(|x^k-\theta_N|\le b_i) -\mathbf{E}_0(\mathbf{I}(|x^k-\theta_N|\le
b_i)\right)|>u/2\}\le A(u)\exp(-B(u)N) \eqno(8)
$$

Denote by $\phi_N(r)=A\exp(-B(r)N)$, where $A$ is a certain constant not depending on $N$. Then it follows from  (6) and (7) that
$$
\Delta(b_i+r)+\phi_N(r)\ge\mathbf{E}_0(\mathbf{I}(|x^k-\theta_N|\le b_i)\ge\Delta(b_i-r)+\phi_N(r)$$

Since the function $\Delta(\cdot)$ satisfies Lipshitz condition (in virtue of boundedness of the density function), from these inequalities for
some $r$ (e.g., $0<r<u/4$) it follows that for large enough $N>N_0$
$$
 \mathbf{P}_0\{|\frac{1}{N}\suml_{k\in\mathcal{N}}\mathbf{I}(|x^k-\theta_N|\le b_i) -\Delta(b_i)|>u/4  \}\le A(u)\exp(-B(u)N)\df\gamma(u,N)
\eqno(9)
$$

Estimate (9) is satisfied for each point $b_i$.

Then in virtue of (4) with the probability no less than $(1-\gamma(u,N))$, for $N>N_0$ we obtain for all $b\in\mathbb{B}_{i}$:
$$
(\Delta(b_i)-u/4)N\le N_1(b)\le (\Delta(b_i)+u/2)N \eqno(10)
$$

Now split the set of all values $N_1(b),b\in\mathbf{B}_{i}$ into two subsets: $\mathbb{A}_i\df\{1\le n\le N: [(\Delta(b_i)-u/4)N]\le n\le
[(\Delta(b_i)+u/2)N]\}$  and its complement. We obtain $\mathbf{P}_0(\mathbf{A}_i)\ge (1-\gamma(u,N))$ при $N>N_0$.

Then
$$\begin{array}{ll}
& \mathbf{P}_0\left\{\supl_{b\in\mathbf{B}_{i}}\{|\suml_{i=1}^{n}\,\tilde x^i|
>\dpfrac C2 n\}\cap \{N_1(b)=n\}\right\}\le \gamma(u,N)\\
&+\mathbf{P}_0\left\{\maxl_{n\in\mathbf{A}_i}\{|\suml_{i=1}^{n}\,\tilde x^i|
>\dpfrac C2 n\}\right\}
\end{array}   \eqno(11)
$$

For the probability in the right hand of (11), we note that $\Delta(b)\ge u$. Hence
$$
\begin{array}{ll}
&\mathbf{P}_0\left\{\maxl_{n\in\mathbf{A}_i}\{|\suml_{i=1}^{n}\,\tilde x^i|
>\dpfrac C2 n\}\right\}\le A(C)\exp\left(-N(\Delta(b_i)-u/4)B(C)\right)\le\\[3mm]
&\le A(C)\exp(-NB(C)3/4u)
\end{array}
\eqno(12)
$$

Since all these considerations are valid for every sub-segment, from (11) and (12) we obtain for each $s=1\,\dots,R$
$$
\mathbf{P}_0\left\{\supl_{b\in\mathbf{B}_s}\{|\suml_{i=1}^{n}\,\tilde x^i|
>\dpfrac C2 n\}\cap \{N_1(b)=n\}\right\}\le A(C)\exp(-NB(C)u/4)
\eqno(13)
$$

The analogous estimate is valid for the second term in (2).

Taking into account (2), (3), (4), (9), (13), we obtain the exponential estimate from theorem 1.

\bigskip
{\bf Theorem 2. Proof.}

Consider the main decision statistic:
$$
\Psi_N(b)=\left(N\suml_{i=1}^{N_1(b)}\,\tilde x^i-N_1(b)\,\suml_{i=1}^N x^i\right)/ N^2.
$$

Write
$$\begin{array}{ll}
& \dpfrac 1N\,\mathbf{E}_{\epsilon}\suml_{i=1}^{N_1}\,\tilde x^i
=\dpfrac 1N\,\suml_{n=1}^N\,\mathbf{E}_{\epsilon}(\suml_{i=1}^n\,\tilde x^i |N_1=n)\,\mathbf{P}_{\epsilon}\{N_1=n\}=\\[3mm]
&=\dpfrac 1N\,\suml_{n=1}^N\,\suml_{i=1}^n\,\mathbf{E}_{\epsilon}(|\tilde x^i-\theta_N|<b | N_1=n)\mathbf{P}_{\epsilon}\{N_1=n\}\\
&=\dpfrac 1N\,\left(\mathbf{E}_{\epsilon} N_1\right)\,\intl_{|\epsilon h-x|<b}\,f(x)xdx / \intl_{|\epsilon h-x|<b}\,f(x)dx \to\intl_{|\epsilon
h-x|<b}\,f(x)xdx, \qquad \text{ при } N\to\infty
\end{array}
$$

Here we used the relationship
$$
\dpfrac 1N \mathbf{E}_{\epsilon}N_1=\dpfrac 1N\,\mathbf{E}_{\epsilon}\suml_{k=1}^{N}k\mathbf{I}(|x^k-\theta_N|\le b)\to \intl_{|\epsilon
h-x|<b}f(x)dx \qquad \text{  при } N\to\infty
$$

Therefore form this relationship from the law of large numbers and the equality
$$
\mathbf{E}_{\epsilon}x^i=\epsilon h
$$
we obtain
$$
\mathbf{E}_{\epsilon}\Psi_N(b)\to \Psi(b) \qquad \text{ при } N\to\infty,
$$
where $\Psi(b)=r(b)-\epsilon h\,d(b)$, $r(b)=\intl_{|\epsilon h-x|<b}f(x)xdx$, $d(b)=\intl_{|\epsilon h-x|<b}f(x)dx$.

For each $C>0$ write:
$$\begin{array}{ll}
\mathbf{P}_{\epsilon}\{|\Psi_N(b)-\Psi(b)| >C\} \le &\mathbf{P}_{\epsilon}\{ |\suml_{i=1}^{N_1(b)}\,\tilde x^i-Nr(b)|
> \dpfrac C2
N\}\\
&+\mathbf{P}_{\epsilon}\{|\dpfrac {N_1(b)}N\,\suml_{i=1}^N\,x^i-N\epsilon h d(b)| > \dpfrac C2 N\}.
\end{array}
\eqno (14)
$$

Consider the first term in the right hand:
$$
\mathbf{P}_{\epsilon}\{ |\suml_{i=1}^{N_1(b)}\,\tilde x^i-Nr(b)|
> \dpfrac C2 N\} \eqno (15)
$$

Write
$$
\mathbf{P}_{\epsilon}\{ |\suml_{i=1}^{N_1(b)}\,\tilde x^i-Nr(b)|
> \dpfrac C2 N\}=P_{\epsilon}\{\suml_{i=1}^{N_1(b)}\,\tilde x^i
<\dpfrac C2 N+Nr(b)\}+P_{\epsilon}\{\suml_{i=1}^{N_1(b)}\,\tilde x^i|<-\dpfrac C2 N+Nr(b)\}.
$$

Consider this estimate for the first probability in the right hand. Write
$$\begin{array}{ll}
&P_{\epsilon}\{\suml_{i=1}^{N_1(b)}\,\tilde x^i <\dpfrac C2 N+Nr(b)\}\le \suml_{n=1}^N\,P_{\epsilon}\{(|\suml_{i=1}^n\,\tilde
x^i|>\dpfrac C2 N+nr(b))\cap (N_1(b)=n)\}\\
&\le \suml_{n=1}^N\, P_{\epsilon}\{|\suml_{i=1}^n\,\tilde x^i|>\dpfrac C2 N+nr(b)\}.
\end{array}
$$

The random value  $\tilde x^i-r(b)$ is centered. Therefore we obtain the following exponential upper estimate:
$$
\suml_{n=1}^N\, P_{\epsilon}\{|\suml_{i=1}^n\,\tilde x^i|>\dpfrac C2 N+nr(b)\}\le L_1\exp(-L_2(C)N),
$$
where the constants $L_1$   and $L_2$ do not depend on $N$.

The second probability in the right hand of (14) is estimated in analogous way. Therefore for any fixed  $b$ we have the following estimate:
$$
P_{\epsilon}\{|\Psi_N(b)-\Psi(b)|>C\}\le L_1 \exp(-L_2(C)N).
$$

The type 2 error probability:
$$\begin{array}{ll}
& P_{\epsilon}\{\maxl_{b\in \mathbf B}|\Psi_N(b)|<C\} \\
& \le P_{\epsilon}\{\maxl_{b\in \mathbf B} |\Psi_N(b)-\Psi(b)|> |\Psi(b)|-C\}.
\end{array}
$$

Let $\delta=\maxl_{b\in \mathbf B}\,|\Psi(b)|-C$. Then
$$
P_{\epsilon}\{\maxl_{b\in \mathbf B}|\Psi_N(b)|<C\}\le L_1\exp( -L_2(\delta) N).
$$

If the sequence of observations satisfies $\psi$-weak dependence condition, we use theorem 1 in order to obtain the exponential estimate
$$
\mathbf{P}_{\epsilon}\{\maxl_{b\in\mathbb{B}}\,|\Psi_N(b)|\le C\}\le L_1\exp(-L_2(\delta)\sqrt{N}).
$$
where $\delta=\maxl_{b\in\mathbb{B}}\,|\Psi(b)|-C
>0$.

As to the proof of 3), remark that the function $\Psi(b)=\mathbf{E}_{\epsilon}\Psi_N(b)$ satisfies the reversed Lipschitz condition in a
neighborhood of $b^*$.

In fact, we have $\Psi(b^*)=0,\,\Psi^{'}(b^*)=0$ and $\Psi^{''}(b^*)=(f(\epsilon h+b^*)-f(\epsilon h-b^*))+b^*(f^{'}(\epsilon
h+b^*)-f^{'}(\epsilon h-b^*))=2(b^*)^2\,f^{''}(u)\ne 0$, where $0\le u=u(b^*)\le b^*$. Therefore in a small neighborhood of $b^*$ we obtain:
$$
|\Psi(b)-\Psi(b^*)|=(b^*)^2\,|f^{''}(u(b^*))|(b-b^*)^2\ge C(b-b^*)^2,
$$
for a certain $C=C(b^*)>0$.

Now for any $0< \kappa <1$ consider the event $|b_N-b^*|>\kappa$. Then
$$
\mathbf{P}_{\epsilon}\{|b_N-b^*|>\kappa\}
 \le \mathbf{P}_{\epsilon}\{\maxl_b\,|\Psi_N(b_N)-\Psi(b^*)| >\dpfrac 12\,C\kappa^2\}\\
\le 4\,\phi_0(\cdot)\,\exp(-L(C)N),
$$
where $L(C)$ is a certain constant not depending on $N$.

 From this inequality it follows that $b_N\to b^*$ $\mathbf{P}_{\epsilon}-$a.s. as $N\to\infty$.

Then
$$
\epsilon_N=N_2(b_N)/N, \qquad h_N=\theta_N/ \epsilon_N
$$
are the nonparametric estimates for $\epsilon$ and $h$, respectively.

In general these estimates are asymptotically biased and non-consistent. For construction of consistent estimates of $\epsilon$ and $h$, we need
information about the d.f. $f_0(\cdot)$. These consistent estimates can be obtained from the following system of equations:
$$\begin{array}{ll}
& \hat\epsilon_N \hat h_N=\theta_N \\
& \dpfrac{1-\hat\epsilon_N}{\hat\epsilon_N}=\dpfrac {f_0(\theta_N-b_N-\hat h_N)-f_0(\theta_N+b_N-\hat
h_N)}{f_0(\theta_N+b_N)-f_0(\theta_N-b_N)}.
\end{array}
$$

The estimates $\hat\epsilon_N$ and $\hat h_N$ are connected with the estimate $b_N$ of the parameter $b^*$ via this system of deterministic
algebraic equations. Therefore the rate of convergence $\hat\epsilon_N \to \epsilon$ and $\hat h_N \to h$ is determined by the rate of
convergence of $b_N$ to $b^*$ (which is exponential w.r.t. $N$). So we conclude that $\hat\epsilon_N \to\epsilon$ and $\hat h_N \to h$
$\mathbf{P}_{\epsilon}$-a.s. as $N\to\infty$.

Theorem 2 is proved.

\end{document}